\documentclass[11pt]{article}
\usepackage{}
\usepackage{enumerate}
\usepackage{mathbbold}
\usepackage{bbm}
\usepackage{amsfonts}
\usepackage[T1]{fontenc}
\usepackage{latexsym,amssymb,amsmath,amsfonts,amsthm}
\usepackage{graphicx, color}
\usepackage{subfigure}
\usepackage{epsf}
\usepackage{epstopdf}
\usepackage{amssymb}
\usepackage{mathrsfs}

\renewcommand{\theequation}{\arabic{section}.\arabic{equation}}

\newcommand{\R}{{\mathbb R}}

\newcommand{\be}{\begin{eqnarray}}
\newcommand{\ben}{\begin{eqnarray*}}
\newcommand{\en}{\end{eqnarray}}
\newcommand{\enn}{\end{eqnarray*}}
\newcommand{\ba}{\backslash}

\newcommand{\ov}{\overline}

\newcommand{\hx}{\hat{x}}
\newtheorem{theorem}{Theorem}[section]
\newtheorem{lemma}[theorem]{Lemma}

\begin{document}
\title{\bf Identification of Electromagnetic Dipoles from Multi-frequency Sparse Electric Far Field Patterns}
\author{Jialei Li\thanks{Academy of Mathematics and Systems Science,
Chinese Academy of Sciences, Beijing 100190, China. Email: lijialei314@buaa.edu.cn}
,\and
Xiaodong Liu\thanks{Academy of Mathematics and Systems Science,
Chinese Academy of Sciences, Beijing 100190, China. Email: xdliu@amt.ac.cn}}
\date{}
\maketitle

\begin{abstract}
The inverse electromagnetic source scattering problem from multi-frequency sparse electric far field patterns is considered. The underlying source is a combination of electric dipoles and magnetic dipoles. We show that the locations and the polarization strengths of the dipoles can be uniquely determined by the multi-frequency electric far field patterns at sparse observation directions. The unique arguments rely on some geometrical discussions and ingenious integrals of the electric far field patterns with properly chosen functions. Motivated by the uniqueness proof, we introduce two indicator functions for locating the magnetic dipoles and the electric dipoles, respectively. Having located all the dipoles, the formulas for computing the corresponding polarization strengths are proposed. Finally, some numerical examples are presented to show the validity and robustness of the proposed algorithm.

\vspace{.2in}

\vspace{.2in} {\bf AMS subject classifications:}
35P25, 45Q05, 78A46, 74B05

\end{abstract}

% ----------------------------------------------------------
% -------------------- New Section -------------------------
% Introduction
\section{Introduction}
\label{sec-intro}
The inverse scattering of acoustic and electromagnetic waves plays an important role in many different areas such as radar, medical imaging, nondestructive testing and geophysical prospection. In this work, we focus on imaging the unknown objects from the measurements taken at sparse sensors. Such a problem arises from many applications where the measurements are difficult or even impossible to be taken all around the unknown objects. Although at a fixed sensor, we can vary the frequency to obtain more data, this is still a small set of data, which brings many difficulties to the problem.
	
The first work in this direction is given by Sylvester and Kelly in \cite{SylvesterKelly}, where they showed a convex polygon containing the unknown source can be uniquely determined from the far field patterns taken at all frequencies, but only at a few observation directions. We refer the reader to \cite{GriesmaierSchmiedecke-source} and \cite{AlaHuLiuSun} for a factorization method and a direct sampling method, respectively, using  multi-frequency far field patterns taken at finitely many observation directions. The direct sampling method has also been applied for locating a perfect conductor with multi-frequency sparse backscattering far field data \cite{ArensJiLiu}. If the multi-frequency sparse scattered fields are measured, the corresponding uniqueness discussion and direct sampling method can be found in a recent work \cite{JiLiu-nearfield}.

Note that with the multi-frequency sparse measurements, the uniqueness of an extended source is still an open problem. If it is known a priori that the object is a combination of point sources, based on some novel geometrical discussions, uniqueness for the locations and the scattering strengths can be established \cite{Ji-multipolarFarfield, JiLiu-point, JiLiu-nearfield}. These uniqueness results also motivates novel direct sampling methods for locating the point sources and the formulas for computing the scattering strengths. A quality-Bayesian approach is proposed in \cite{LiLiuSunXu, LiuGuoSun} to reconstruct the locations and intensities of the unknown acoustic sources.  We also refer to \cite{GriesmaierSchmiedecke} for an MUSIC-type reconstruction method for small inhomogeneities
by multi-frequency sparse far field patterns.

In this paper we consider the inverse source scattering problem for time harmonic electromagnetic waves using multi-frequency sparse electric far field patterns.
For an extended source, it is shown in \cite{JiLiu-electromagnetic} that the smallest strip containing the source support with the observation direction as the normal can be uniquely determined by the multi-frequency far field pattern at a single observation direction. Similar to the acoustic case, a convex support of the extended source is then expected to be reconstructed with the data in finitely many observation directions. This paper focus on the case that the source to be determined is a combination of electric dipoles and magnetic dipoles. Different to the general source functions, the dipoles are characterized by their locations and polarization strengths. This enables us to study the uniqueness analyses and numerical algorithms with less data. This is also a nontrivial extension of the results  \cite{Ji-multipolarFarfield, JiLiu-point} for inverse acoustic source scattering problems. Actually, in contrast to the acoustic case, the electric field generated by different types of dipoles are coupled together, which makes the problem more complicated. Inspired by \cite{Ji-multipolarFarfield}, to distinguish the dipoles, we introduce two weighted integral using electric far field pattern at two opposite observation directions. Another difficulty is due to the fact that the electric far field pattern vanishes if the observation direction is parallel to the polarization direction of the dipole. Hence, more properly chosen directions are necessary to determine the dipoles when the polarization directions are unknown. Note that the electromagnetic fields generated by the electric dipoles and magnetic dipoles can be viewed as the fundamental solution of the Maxwell equations \cite{CK}, thus the results obtained in this paper are expected to be helpful for understanding the inverse scattering problems with extended sources or scatterers.
	
The remaining part of this paper is organized as follows. In the next section, we fix the notations and introduce the scattering of the mixed electric dipoles and magnetic dipoles. Some uniqueness results are shown in section \ref{sec-Uniq-ElecMag}. Motivated by the uniqueness arguments, Section \ref{sec-Algorithm} is devoted to introduce two indicator functions for locating the two types of dipoles, respectively, and the numerical formulas for computing the corresponding polarization strengths.  Numerical examples are presented in section \ref{sec-NumericalExample} to verify the proposed numerical algorithm.

% ----------------------------------------------------------
% -------------------- New Section -------------------------
% The Direct Electromagnetic point sources problem
\section{Scattering due to electromagnetic dipoles}
\label{sec-directProb}
We consider the scattering of the electromagnetic point sources in a homogeneous isotropic medium in $\R^3$. The point sources are a combination of magnetic dipoles and electric dipoles.
	
% subsec- electric field of magnetic dipoles
% \subsection{Electromagnetic field of a magnetic dipole}
% \label{subsec-ElecMagField-AMag}

Let $k>0$ be the wave number.
For a magnetic dipole located at $z\in \mathbb{R}^3$, the corresponding electromagnetic field is given by \cite{CK}
\be\label{MagneticDipole}
E_{mag}(x,z,q,k) :=  {\rm curl}_x\ q \Phi_k (x,z),
\quad
H_{mag}(x,z,q,k) :=  \frac{1}{ik}{\rm curl}_x E_{mag}(x,z,q,k),
\quad x\neq z,
\en
where
\ben\label{FundaSolu_of_3DimHelmholtz}
\Phi_k(x,z)=\frac{e^{ik|x-z|}}{4\pi |x-z|},\,\ x\neq z,
\enn
is the fundamental solution to the Helmholtz equation $\Delta u +k^2u=0$, $q :=\tau p$ is the polarization strength given by product of the strength $\tau\in \mathbb{C}$ and the polarization direction $p\in S^2:= \{ x\in \R^3|\quad |x|=1\}$.
Similarly,
\be\label{ElectricDipole}
H_{elec}(x,z,q,k) :=  {\rm curl}_x\ q \Phi_k (x,z),
\quad
E_{elec}(x,z,q,k) :=  \frac{i}{k}{\rm curl}_x H_{elec}(x,z,q,k),\quad x\neq z
\en
represent the electromagnetic field generated by an electric dipole.	
%In terms of \eqref{MagneticDipole}-\eqref{ElectricDipole}, the dipole $(\tau,p,z)$ and the dipole $(-\tau,-p,z)$ generate exactly the same electric field. Hence we will always use $q := \tau p$ rather than $\tau p$ without extra explanation. In fact, there is a one-to-one correspondence between $q=\tau p$ and $\{(\tau, p), (-\tau,-p)\}$. Therefore we say the dipole is uniquely determined if $z_m$ and $q_m$ are uniquely determined.
As mentioned in page 230 of \cite{CK}, the fields \eqref{MagneticDipole} and \eqref{ElectricDipole} can be viewed as the fundamental solution of the Maxwell equations
\ben
{\rm curl} E-ik H=0,\quad {\rm curl} H+ikE = 0.
\enn

By straightforward calculations, it can be seen that
\ben
E_{mag}(x,z,q,k) = \frac{e^{ik|x|}}{|x|} \left\{ E_{mag}^\infty (\hat{x},z,q,k) + O\left( \frac{|q|}{|x|} \right) \right\},\quad x\in\R^3\ba\{z\},
\enn
uniformly in all directions $\hx=x/|x|$	with
\be
\label{ElecFarField_of_SingleMagDipole}
E_{mag}^\infty (\hat{x},z,q,k)=  \frac{ik}{4\pi} e^{-ik \hat{x}\cdot z}\hat{x}\times q, \quad \hat{x}\in S^2,
\en
and
\ben
E_{elec}(x,z,q,k) = \frac{e^{ik|x|}}{|x|} \left\{ E_{elec}^\infty (\hx,z,q,k) + O\left( \frac{|q|}{|x|} \right) \right\},\quad x\in\R^3\ba\{z\},
\enn
uniformly in all directions $\hx=x/|x|$	with
\be\label{ElecFarField_of_SingleElecDipole}
E_{elec}^\infty (\hx,z,q,k) =  \frac{ik}{4\pi} e^{-ik \hat{x}\cdot z} \hat{x}\times (q\times \hat{x}), \quad \hat{x}\in S^2.
\en
Here, $E_{mag}^\infty$ and $E_{elec}^\infty$ defined on the unit sphere $S^2$ are known as the electric far field patterns for the magnetic dipole and the electric dipole, respectively.

We consider an array of $M$ dipoles at locations $z_1, z_2, \cdots, z_M\in\R^3$ in  the homogeneous isotropic space $\R^3$. Denote by $q_m := \tau_m p_m \in \mathbb{C}^3$ the polarization strength vector of the $m$-th dipole.
Without loss of generality we assume that
the first $M_1$ dipoles are magnetic dipoles and the other $M_2 := M-M_1$ dipoles are electric dipoles.
Then, the electric field due to $M$ dipoles is given by
\ben
&&E (x,k)\cr
&=& \sum_{m=1}^{M_1} E_{mag}(x,z_m,q_m,k) + \sum_{m=M_1+1}^M E_{elec}(x,z_m,q_m,k) \cr
&=& \sum_{m=1}^{M_1} {\rm curl}_x\ q_m \Phi_k (x,z_m) + \sum_{m=M_1+1}^M \frac{i}{k}{\rm curl}_x{\rm curl}_x\ q_m \Phi_k (x,z_m),\quad x\in\R^3\ba\{z_1, z_2, \cdots, z_M\}.
\enn
Correspondingly, the electric far field pattern is given by
\be
\label{eq-FarFOfElecAndMagDipoles}
E^\infty (\hat{x},k) =
\sum_{m=1}^{M_1} \frac{ik}{4\pi} e^{-ik \hat{x}\cdot z_m} \hat{x} \times q_m +
\sum_{m=M_1+1}^M \frac{ik}{4\pi} e^{-ik \hat{x}\cdot z_m} \hat{x} \times (q_m \times \hat{x}),\quad \hx\in S^2.
\en

Denote by
\be\label{ThetaL}
\Theta_L=\{\hat{x}_l \in S^2\, |\, l = 1,2,\cdots, L\}
\en
the set of $L$ properly chosen observation directions. The {\em inverse scattering problem} is to identify the dipoles from the multi-frequency electric far field patterns $E^\infty (\hat{x},k)$
at sparse observation directions $\hx\in\Theta_L$. More specifically, we want to determine the locations $z_1, z_2, \cdots, z_M$, to reconstruct the polarization strength vectors $q_1, q_2, \cdots, q_M$ and to clarify the dipole type of each dipole.

% ----------------------------------------------------------
% -------------------- New Section -------------------------
% Uniqueness of Inverse Electromagnetic Point Sources Problem
\section{Uniqueness of inverse electromagnetic point sources problem}
\label{sec-Uniq-ElecMag}
In this section, we investigate under what conditions the dipoles can be identified by the multi-frequency sparse electric far field patterns.

% case $M=M_1=1$
\subsection{Uniqueness when $M=M_1=1$}
\label{subsec-Unique-AMag}
	We begin with the simplest case that identifying a single magnetic dipole, i.e., $M=M_1=1$.
%By \eqref{MagneticDipole} ,  the electric scattered field is
%\be
%\label{eq-SingleMag-Near}
%E(x,k) &=  \frac{ik|x-z_1|-1}{4\pi |x-z_1|^2} e^{ik|x-z_1|}\
%			\frac{x\times q_1}{|x-z_1|}, \quad x\in \mathbb{R}^3, x\neq z.
%\en

% Theorem $M=M_1=1$ Uniqueness fixed k
\begin{theorem}\label{Unique-SingleMag}
For a fixed wave number $k>0$, when $M=M_1=1$, we have the following results:
\begin{itemize}
\item[(1)] Given the location $z_1$, the polarization strength $q_1$ can be uniquely determined by $E^\infty(\hat{x}_1,k),E^\infty(\hat{x}_2,k)$ where $\hat{x}_1, \hat{x}_2$ are two non-collinear observation directions in $S^2$;
%% ---
%\item[(b)] Given the location $z_1$, let $x_1, x_2\in \mathbb{R}^3$ such that $x_1 \times x_2 \neq 0$, then the polarization strength $q_1$ can be uniquely determined by $E(x_1, k), E(x_2, k)$;
%% ---
%\item[(c)] Given the polarization strength $q_1$, for fixed $x\in \mathbb{R}^3\setminus \{z_1\}$ such that $x \times q_1 \neq 0$, the distance $r := |x-z_1|$ of $x$ and $z_1$ can be uniquely determined by $|E(x,k)|$. Moreover, $z_1$ can be uniquely determined by $|E(x_j,k)|, j=1,2,3,4$, where $\{x_j\neq z_1|x_j \times q_1 \neq 0, j=1,2,3,4\}$ are four sensors that are not coplanar;
%% ---
\item[(2)] Given the polarization strength $q_1$, suppose that $|k \hat{x} \cdot z_1|<\pi$ and $\hat{x} \times q_1\neq 0$, then $\hat{x} \cdot z_1$ can be uniquely determined by $E^\infty(\hat{x}, k)$. Moreover, $z_1$ can be uniquely determined by $E^\infty(\hat{x}_j, k)$ for three linearly independent observation directions $\hat{x}_j, j=1,2,3$.
\end{itemize}
\end{theorem}
% begin of proof
\begin{proof}
\begin{itemize}
\item[(1)]
   %The vectors $\hx$, $q_1\times \hx$ and $\hx\times(q_1\times\hx)$ form a basis in $\R^3$ provided $q_1\times \hx\neq 0$. Noting that $\hat{x}_1\neq \pm \hat{x}_2$, we have $\hx_1$ and $\hx_2$ are not collinear and without of loss of generality we assume that  $q_1\times \hx_1\neq 0$.
   With the help of the vector identity $a\times(b\times c)= (a\cdot c) b-(a\cdot b)c$, noting the fact that $\hx_1\in S^2$, we have
   \be\label{q1}
   q_1 = (q_1\cdot\hx_1)\hx_1+\hx_1\times(q_1\times\hx_1).
   \en
   Taking the vector product of \eqref{q1} with $\hx_2$ yields that
   \ben
   \hx_2\times q_1 = (q_1\cdot\hx_1)\hx_2\times\hx_1 +\hx_2\times[\hx_1\times(q_1\times\hx_1)].
   \enn
Noting that $\hat{x}_1$ and $\hat{x}_2$ are non-collinear, we have $\hx_1\times\hx_2\neq 0$. Therefore,
\be\label{q1dothx1}
q_1\cdot\hx_1 = \frac{\hat{x}_2\times \hat{x}_1}{|\hat{x}_2\times \hat{x}_1|^2} \cdot \left(
\hat{x}_2 \times q_1 - \hx_2\times[\hx_1\times(q_1\times\hx_1)]\right).
\en
With the help of \eqref{ElecFarField_of_SingleMagDipole}, we have
   \be\label{hxtimesq}
   \hat{x}_j \times q_1=\frac{4\pi}{ik} e^{ik \hat{x}_j\cdot z_1}E^\infty(\hat{x}_j,k),\quad j=1,2.
   \en
Inserting \eqref{q1dothx1} and \eqref{hxtimesq} into \eqref{q1} we derive that
\be\label{q1formula}
q_1 &=& \frac{4\pi\hat{x}_2\times \hat{x}_1}{ik|\hat{x}_2\times \hat{x}_1|^2} \cdot \left(
e^{ik \hx_2\cdot z_1}E^\infty(\hx_2,k) + e^{ik \hat{x}_1\cdot z_1}\hx_2\times[\hx_1\times E^\infty(\hx_1,k)]\right)\hx_1\cr
&&-\frac{4\pi}{ik} e^{ik \hat{x}_1\cdot z_1}\hx_1\times E^\infty(\hx_1,k),
\en
which is actually the formula for computing $q_1$ from the electric far field patterns $E^{\infty}(x_j,k), \,j=1,2$.
%% ---
%\item[(b)] This follows from similar arguments as in the proof of (a) by replacing \eqref{hxtimesq} by
%\ben
% x_j \times q_1 = \frac{4\pi |x_j - z_1|^3}{ik|x-z_1| -1} e^{-ik|x-z_1|} E(x,k), j=1,2.
%\enn
%% ---
%\item[(c)] Note that the right hand side of
%\ben
%\left | E(x,k)\right|^2 = \frac{1+k^2 r^2}{16\pi^2 r^6} |x\times q_1|^2
%\enn
%is strictly decreasing with respect to $r>0$. Therefore $r=|x-z_1|$ is uniquely determined from $|E(x,k)|$ at the fixed sensor $x\in\R^3\ba\{z_1\}$.
%By a geometrical discussion (see e.g., the proof of Theorem 3.1 in \cite{JiLiu-point}),
%we have that $z_1$ is uniquely determined by $|E(x_j,k)|$ with four non-coplanar sensors $x_j, j=1,2,3,4$.
%% ---
\item[(2)] With the help of the representation \eqref{ElecFarField_of_SingleMagDipole}, we have
\ben
e^{-ik\hat{x}\cdot z_1}=\frac{4\pi}{ik} E^\infty(\hat{x},k)\cdot \frac{\hat{x}\times q_1}{|\hat{x}\times q_1|^2},\quad \hx\in S^2.
\enn
This implies that $\hat{x}\cdot z_1$ is uniquely determined by $E^\infty(\hat{x},k)$ when $|k\hat{x}\cdot z_1|<\pi$. The uniqueness of $z_1$ then follows by noting that the three observation directions $\hat{x}_j, j=1,2,3$ are linearly independent.
\end{itemize}
\end{proof}
% end of proof

Note that, in the above local uniqueness results, we have only used the electric far field pattern for a fixed wave number $k$. The following theorem shows that both the location $z_1$ and polarization strength $q_1$ can be uniquely determined when multiple frequencies are used.

% theorem-uniqueness multi-frequency
\begin{theorem}\label{Unique-SingleMag-multiK}
Let $M=M_1 = 1$, the wave numbers be contained in some interval $0< k_-<k<k_+<+\infty$ and $\hat{x}_1, \hat{x}_2$ and $\hat{x}_3$ be three linearly independent observation directions such that $\frac{1}{k_-} E^\infty(\hat{x}_j, k_-)\neq \frac{1}{k_+} E^\infty(\hat{x}_j,k_+)$ for $j=1,2,3$. Then both the dipole location $z_1$ and its polarization strength $q_1$ can be uniquely determined by $E^\infty (\hat{x}_j,k), k\in (k_-,k_+), j=1,2,3$.
\end{theorem}
% begin of proof
\begin{proof}
Under the assumption that $\frac{1}{k_-} E^\infty(\hat{x}, k_-)\neq \frac{1}{k_+} E^\infty(\hat{x},k_+)$ we have
\be\label{hxq1z1neq0}
\hx\times q_1\neq 0\quad\mbox{and}\quad \hx\cdot z_1\neq 0\quad\mbox{for}\quad \hx\in\{\hat{x}_1, \hat{x}_2, \hat{x}_3\}.
\en
Using the far field representation \eqref{ElecFarField_of_SingleMagDipole}, we have
\be\label{eq-SingleMag-findExp}
\frac{1}{k}E^\infty(\hat{x},k) \cdot \overline{E^\infty}(\hat{x},k_-)
= \frac{k_-}{16 \pi^2} e^{-i(k-k_-) \hat{x} \cdot z_1}
|\hat{x}\times q_1|^2, \quad \hx\in S^2, k\in (k_-, k_+),
\en
where $\overline{\cdot}$ denotes the complex conjugate.
Integrating the above equation with respect to $k$ over $(k_-, k_+)$, in terms of \eqref{hxq1z1neq0}, we obtain
\ben
&&\int_{k_-}^{k_+} \frac{1}{k} E^\infty(\hat{x},k) \cdot \overline{E^\infty}(\hat{x},k_-) {\rm d} k\cr
&=& \frac{ k_- |\hat{x}\times q_1|^2}{16 \pi^2}\int_{k_-}^{k_+} e^{-i(k-k_-) \hat{x} \cdot z_1}{\rm d} k \cr
&=& \frac{1}{-i \hat{x}\cdot z_1}\frac{ k_- |\hat{x}\times q_1|^2}{16 \pi^2}
\left[e^{-i(k_+-k_-) \hat{x} \cdot z_1} -1\right] \cr
&=& \frac{1}{-i \hat{x}\cdot z_1}
\left[
\frac{1}{k_+} E^\infty(\hat{x},k_+) \cdot \overline{E^\infty}(\hat{x},k_-)
- \frac{1}{k_-} E^\infty(\hat{x},k_-) \cdot \overline{E^\infty}(\hat{x},k_-)
\right],\quad \hx\in\{\hat{x}_1, \hat{x}_2, \hat{x}_3\}.
\enn
Note also that $\int_{k_-}^{k_+} \frac{1}{k} E^\infty(\hat{x},k) \cdot \overline{E^\infty}(\hat{x},k_-) {\rm d}k \neq 0$ by the inequality \eqref{hxq1z1neq0}. Therefore
\be
\label{eq-SingleMag-hatx-dot-z}
\hat{x}\cdot z_1
= i \frac{k_- E^\infty(\hat{x},k_+) \cdot \overline{E^\infty}(\hat{x},k_-)- k_+ E^\infty(\hat{x},k_-) \cdot \overline{E^\infty}(\hat{x},k_-)}
{k_+k_-\int_{k_-}^{k_+} \frac{1}{k} E^\infty(\hat{x},k) \cdot \overline{E^\infty}(\hat{x},k_-) {\rm d}k},\quad \hx\in\{\hat{x}_1, \hat{x}_2, \hat{x}_3\},
\en
is uniquely determined from $E^\infty (\hat{x},k), k\in (k_-,k_+)$. Furthermore, $z_1$ is uniquely determined from $E^\infty (\hat{x}_j,k), k\in (k_-,k_+), j=1,2,3$ by the fact that $\hat{x}_1, \hat{x}_2$ and $\hat{x}_3$ are three linearly independent directions. Finally, the uniqueness of $q_1$ is obtained by the first result in Theorem \ref{Unique-SingleMag}.
\end{proof}
% end of proof

After slight modifications of the proof, the analogous uniqueness results of Theorem \ref{Unique-SingleMag} and Theorem \ref{Unique-SingleMag-multiK} can be formulated for a single electric dipole. An interesting question is how to distinguish the dipole types from the electric far field patterns. We will give an answer in the next subsection for more complex multiple mixed dipoles.
	
% Many Dipoles
\subsection{Uniqueness for multiple mixed dipoles}\label{subsec-Unique-ManyDipoles}
In this subsection we study the much more complicated case with multiple mixed dipoles.
As shown in \eqref{ElecFarField_of_SingleMagDipole}-\eqref{eq-FarFOfElecAndMagDipoles}, an obvious difficulty is how to decouple the electric dipoles and the magnetic dipoles. Besides, for multiple dipoles, the proper choice of the observation directions plays an important role in locating the dipoles and clarifying the dipole type.
%	It's much more complicated when there are more dipoles because there are two types of dipoles. Moreover, the nonlinearity between the measurement and the location of dipoles brings difficulties. In the following theorem, we shows that multi-frequency electric far field pattern data uniquely determined the dipoles with enough observation directions.

Recall the observation direction set $\Theta_L$ given in \eqref{ThetaL}. Define by
\be
\label{Definition-Planes}
\Pi_{l,m}= \{ z\in \mathbb{R}^3\,|\, \hat{x}_l \cdot (z-z_m)=0\},\quad
1\leq l \leq L, \quad 1\leq m \leq M,
\en
the plane passing through $z_m$ with normal $\hat{x}_l \in \Theta_L$. For any point $z\in \mathbb{R}^3$,
denote by
\ben
f(z):= {\mathrm Card} \{\Pi_{l,m}\,|\, z\in \Pi_{l,m},\, 1\leq l\leq L,\, 1\leq m\leq M\}
\enn
the number of the planes $\Pi_{l,m}, 1\leq l\leq L, 1\leq m\leq M$ passing through $z$.
%, where ${\mathrm Card}$ is the cardinal of the set, i.e., the number of the elements (for finite set).

% lemma-upper bound of $f(z)$ when $z$ is not the location of any dipole
\begin{lemma}
\label{lemma-UpperBound-f(z)}
Giving $M$ points $z_1, z_2, \cdots, z_M$ in $\mathbb{R}^3$. Assume that any three directions in $\Theta_L$ are not coplanar, then
\ben
f(z)
\left\{
  \begin{array}{ll}
    \leq 2M, & \hbox{$z\in \mathbb{R}^3 \setminus \{z_1,z_2, \cdots, z_M\}$;} \\
    = L, & \hbox{$z\in \{z_1,z_2, \cdots, z_M\}$.}
  \end{array}
\right.
\enn
\end{lemma}
% begin of proof
\begin{proof}
$f(z_m) = L , 1\leq m\leq M$ follow immediately from the fact that $z_m \in \Pi_{l,m}$ for all $1\leq l\leq L$.
	
Now we proof $f(z)\leq 2M$ for $z\in \mathbb{R}^3 \setminus \{z_1,z_2, \cdots, z_M\}$. Assume to the contrary that $f(z^*)\geq 2M+1$ for some $z^* \in \mathbb{R}^3 \setminus \{z_1,z_2, \cdots, z_M\}$. Then, by the pigeonhole principle, there exists a point $z_m$ such that at least three planes $\Pi_{l_1,m},\Pi_{l_2,m},\Pi_{l_3,m}$ pass through $z^*$, i.e.,
\ben
\hat{x}_{l_j} \cdot(z^*-z_m) = 0, \quad j =1, 2, 3.
\enn
However, $\hat{x}_{l_j}, j=1,2,3$ are not coplanar, which implies $z^*-z_m=0$. This is a contradiction to $z^* \in \mathbb{R}^3 \setminus \{z_1,z_2, \cdots, z_M\}$.
%Therefore $f(z)\leq 2M$ for all $z\in\mathbb{R}^3 \setminus \{z_1,z_2, \cdots, z_M\}$.
\end{proof}
% end of proof

Lemma \ref{lemma-UpperBound-f(z)} shows that dipoles can be uniquely determined when $L\geq 2M+1$ if $f(z)$ can be calculated for all $z\in \mathbb{R}^3$ from the electric far field patterns. By \eqref{eq-FarFOfElecAndMagDipoles}, for dipole with polarization direction $p_m$, the electric far field at the observation direction $\pm p_m$ is zero, which means they contain no information about the dipole. Therefore, given dipole $(z_m, \tau_m, p_m)$ and observation $\hat{x}_l$, it might be impossible to calculate $\hat{x}_l\cdot z_m$ and to obtain $\Pi_{l,m}$ from the electric far field pattern.
Fortunately, the number of observation directions that don't work is finite, which indicates uniqueness is possible when the number of observation directions $L$ is large enough.

For $z\in\R^3$, we define by
\be
\label{Definition-CountFunc}
\begin{aligned}
f_{mag}(z) &:= {\mathrm Card} \{\Pi_{l,m}\,|\, z\in \Pi_{l,m},\, 1\leq l\leq L,\, 1\leq m\leq M_1\} \\
f_{elec}(z) &:= {\mathrm Card} \{\Pi_{l,m}\,|\, z\in \Pi_{l,m},\, 1\leq l\leq L,\, M_1< m\leq M\}
\end{aligned}
\en
the numbers of the planes $\{\Pi_{l,m}| 1\leq m \leq M_1\}$ and $\{\Pi_{l,m}| M_1< m \leq M\}$ that passing through $z$, respectively.	
% theorem-uniqueness-manydipoles-multi-frequency
\begin{theorem}
\label{Theorem-ElecMagUnique-multi-frequency}
If any three directions in $\Theta_L$ are not coplanar and $L\geq \max\{4M_1, 4(M-M_1)\}$, then
\ben
\{(z_m,q_m)\,|\, 1\leq m \leq M_1\}\quad\mbox{and}\quad \{(z_m,q_m)\,|\, M_1< m \leq M\}
\enn
are uniquely determined by multi-frequency electric far field patterns $E^\infty(\pm \hat{x},k), \hat{x}\in \Theta_L, k\in (0,\infty)$.
\end{theorem}
% begin of proof
\begin{proof}
Given $\hat{x}\in \Theta_L$ and $z\in\R^3$, define
\ben
\begin{aligned}
f_{\hat{x},mag}(z) &= \left\{
\begin{aligned}
&1,  &\mbox{if } \ \hat{x}\cdot (z-z_m)=0  \mbox{ for some } m\in \{1,2,\cdots, M_1\}; \\
&0,  &\mbox{otherwise},
\end{aligned}
\right.\\
f_{\hat{x},elec}(z) &= \left\{
\begin{aligned}
&1,  &\mbox{if } \ \hat{x}\cdot (z-z_m)=0  \mbox{ for some } m \in \{M_1+1,M_1+2,\cdots, M\}; \\
&0,  &\mbox{otherwise}.
\end{aligned}
\right.
\end{aligned}
\enn
We have
\ben
f_{mag}(z) = \sum_{\hat{x}\in \Theta_L} f_{\hat{x},mag}(z)\quad\mbox{and}\quad f_{elec}(z) = \sum_{\hat{x}\in \Theta_L} f_{\hat{x},elec} (z)\quad\mbox{for}\quad z\in\R^3.
\enn

For any $K>0$, in view of \eqref{eq-FarFOfElecAndMagDipoles} we have
\be\label{eq-integralOfMag}
&&F_{mag}(z,\hat{x}, K)\cr
&:=& \frac{2\pi}{K}\int_0^K \frac{1}{ik}\left[e^{ik\hat{x}\cdot z} E^\infty(\hat{x},k) - e^{-ik\hat{x}\cdot z} E^\infty(-\hat{x},k)\right]{\rm d} k \cr
&=& \sum_{m=1}^{M_1}  \frac{1}{2K} \int_{-K}^K  e^{ik \hat{x}\cdot (z-z_m)}  {\rm d} k\  \hat{x} \times q_m \cr
&&\quad +\sum_{m=M_1+1}^M  \frac{1}{2K} \int_{0}^K \left(e^{ik \hat{x}\cdot (z-z_m)}-e^{-ik \hat{x}\cdot (z-z_m)}\right)  {\rm d}{k}\  \hat{x} \times (q_m \times\hat{x}) \cr
&=& \sum_{m=1}^{M_1}  \frac{1}{2K} \int_{-K}^K  e^{ik \hat{x}\cdot (z-z_m)}  {\rm d}{k}\  \hat{x} \times q_m \cr
&&\quad +\sum_{m=M_1+1}^M  \frac{i}{K} \int_{0}^K \sin(k \hat{x}\cdot(z-z_m))
  {\rm d}{k}\  \hat{x} \times (q_m \times \hat{x}),\quad z\in \mathbb{R}^3, \hat{x}\in \Theta_L.\quad
\en
Letting $K\rightarrow\infty$, we see that
\be
\label{eq-limit-of-Fmag}
F_{mag}(z,\hat{x}) :=\lim_{K\to \infty} F_{mag}(z,\hat{x}, K) =
\sum_{\hat{x}\cdot(z-z_m)=0,\,  1\leq m\leq M_1} \hat{x}\times q_m, \quad z\in \mathbb{R}^3, \hat{x}\in \Theta_L.
\en
%If $F_{mag}(z,\hat{x})\neq 0$, $\hat{x}\cdot (z-z_m)=0$ for some $m\leq M_1$, which implies $f_{\hat{x}, mag}(z)=1$.
Define
\be\label{eq-Definition-tilde-f-mag-singledirection}
\tilde{f}_{\hat{x},mag} (z) :=  \left\{
\begin{aligned}
&1,  &\mbox{if }F_{mag}(z,\hat{x})\neq 0; \\
% ------
&0,  &\mbox{otherwise}.
\end{aligned}
\right.
\en
Obviously, $f_{\hat{x},mag}(z)\geq \tilde{f}_{\hat{x},mag}(z)$. By Lemma \ref{lemma-UpperBound-f(z)}, we have
\be\label{zothers}
\sum_{\hat{x}\in \Theta_L} \tilde{f}_{\hat{x},mag}(z)\leq \sum_{\hat{x}\in \Theta_L} f_{\hat{x},mag}(z) = f_{mag}(z)\leq 2M_1\leq L/2,\quad z\in \mathbb{R}^3\setminus \{z_1,z_2, \cdots, z_{M_1}\}.
\en
On the other hand, for fixed $z_m, 1\leq m\leq M_1$, there are at most $2(M_1-1)$ observation directions $\hat{x}$ satisfying $\hat{x}\cdot (z_m - z_{m^*})=0$ for some $m^*\neq m, 1\leq m^*\leq M_1$ and at most one direction $\hat{x}$ satisfying $\hat{x}\times q_m = 0$. Hence we have at least $L- 2M_1+1$ observation directions $\hat{x}$ such that
\ben
F_{mag}(z_m,\hat{x}) = \hat{x} \times q_m \neq 0,
\enn
which implies $\tilde{f}_{\hat{x},mag}(z_m)=1$. Consequently we have
\be\label{zmag}
\sum_{\hat{x}\in \Theta_L} \tilde{f}_{\hat{x},mag}(z)\geq L-2M_1+1 > L/2,\quad z\in \{z_1,z_2, \cdots, z_{M_1}\}.
\en
In terms of \eqref{zothers} and \eqref{zmag}, we deduce that $\{z_m| 1\leq m\leq M_1\}=\{z\in \mathbb{R}^3|\sum_{\hat{x}\in \Theta_L} \tilde{f}_{\hat{x},mag}(z)>L/2\} $ is uniquely determined.

For each magnetic dipole $(z_m,q_m),\, 1\leq m \leq M_1$, we can always find at least two linearly independent observation directions $\hx_1$ and $\hx_2$ such that \ben
F_{mag}(z_m,\hx_j) = \hx_j \times q_m, \quad j=1,2.
\enn
Following the arguments in the proof of Theorem \ref{Unique-SingleMag} (1) we have
\be\label{qmmagformula}
q_m &=& \frac{\hat{x}_2\times \hat{x}_1}{|\hat{x}_2\times \hat{x}_1|^2} \cdot \Big(
F_{mag}(z_m,\hx_2) + \hx_2\times[\hx_1\times F_{mag}(z_m,\hx_1)]\Big)\hx_1\cr
&&-\hx_1\times F_{mag}(z_m,\hx_1),\quad 1\leq m \leq M_1.
\en
Therefore the magnetic dipoles $\{(z_j,q_j)| 1\leq j \leq M_1\}$ are uniquely determined.

For electric dipoles, we define
\be
F_{elec}(z,\hat{x}) &:=&\lim_{K\to \infty} F_{elec}(z,\hat{x},K),\quad z\in \mathbb{R}^3, \hat{x}\in \Theta_L
\en
with
\be
\label{IntegralOfElec}
F_{elec}(z,\hat{x},K) &:=& \frac{2 \pi}{K} \int_0^K \frac{1}{ik}\left[
e^{ik\hat{x}\cdot z} E^\infty(\hat{x},k) + e^{-ik\hat{x}\cdot z} E^\infty(-\hat{x},k)
\right]{\rm d}{k} ,\quad z\in \mathbb{R}^3, \hat{x}\in \Theta_L. \qquad
\en
%and
%\ben
%\tilde{f}_{\hat{x},elec} (z) :=  \left\{
%\begin{aligned}
%&1,  &\mbox{if }F_{elec}(z,\hat{x})\neq 0; \\
%% -----
%&0,  &\mbox{otherwise}.
%\end{aligned}
%\right.
%\enn
Inserting \eqref{eq-FarFOfElecAndMagDipoles} into \eqref{IntegralOfElec} and letting $K\rightarrow\infty$ we have
\be
\label{eq-limit-of-Felec}
F_{elec}(z,\hat{x}) = \sum_{\hat{x}\cdot(z-z_m)=0,\, M_1<m\leq M} \hat{x}\times (q_m \times \hat{x}),\quad z\in \mathbb{R}^3, \hat{x}\in \Theta_L.
\en
The uniqueness of the locations $\{z_{M_1+1}, z_{M_1+2},\cdots, z_M\}$ for the electric dipoles follows from similar arguments for the magnetic dipoles.

For each $z_m\in\{z_{M_1+1}, z_{M_1+2},\cdots, z_M\}$, similar to the case of magnetic dipoles, we have at least two linearly independent observation directions $\hat{x}_1$ and $\hat{x}_2$ such that
\be\label{Feleczm}
F_{elec}(z_m,\hat{x}) = \hat{x}\times (q_m \times \hat{x}),\quad \hat{x}\in\{\hat{x}_1, \hat{x}_2\}.
\en
In view of \eqref{q1} we have
\be\label{qmelec}
   q_m=(q_m\cdot\hx_1)\hx_1+F_{elec}(z_m,\hx_1).
\en
From this and \eqref{Feleczm} we obtain that
   \ben
   \hx_2\times F_{elec}(z_m,\hx_2) = \hx_2\times q_m = (q_m\cdot\hx_1)(\hx_2\times \hx_1) +\hx_2\times F_{elec}(z_m,\hx_1).
   \enn
Therefore,
\ben
q_m\cdot\hx_1 = \frac{\hx_2\times \hx_1}{|\hx_2\times \hx_1|^2} \cdot \Big(\hx_2\times F_{elec}(z_m,\hx_2) - \hx_2\times F_{elec}(z_m,\hx_1)\Big).
\enn
Inserting this into \eqref{qmelec}, we derive that
\be\label{qmelecformula}
q_m = \frac{\hx_2\times \hx_1}{|\hx_2\times \hx_1|^2} \cdot \Big(\hx_2\times F_{elec}(z_m,\hx_2) - \hx_2\times F_{elec}(z_m,\hx_1)\Big)\hx_1+F_{elec}(z_m,\hx_1),
\en
which is a formula for computing $q_m$ from $F_{elec}(z_m,\hx_1)$ and $F_{elec}(z_m,\hx_2)$. This completes the proof of the theorem.
\end{proof}

Finally, we want to remark that less data is needed provided a priori information on the physically property of the dipoles. For example, assume that all the point sources $(z_m, q_m)$ are magnetic dipoles with $q_m\in\R^3,\,1\leq m \leq M$. Then
\ben
E^{\infty}(-\hx,k) = \ov{E^{\infty}(\hx,k)},\quad \hx\in S^2.
\enn
Therefore, the magnetic dipoles are uniquely determined from only half of the data used in Theorem \ref{Theorem-ElecMagUnique-multi-frequency}, i.e, the multi-frequency electric far field patterns $E^\infty(\hat{x},k), \hat{x}\in \Theta_L, k\in (0,\infty)$.

% special case
\subsection{Uniqueness when dipoles are in a plane}
\label{subsec-Unique-ElecMagInAPlane}
In this subsection, we consider an interesting case that all of the dipoles are located in a known plane. Less data is then needed to identify the dipoles. Without lose of generality, suppose that the known plane is
\ben
\Pi:=\{z=(z^1, z^2, z^3)\in \mathbb{R}^3\, |\, z^3 = 0\}.
\enn
%(otherwise we change the coordinates).
It's interesting that the observation directions set $\Theta_L$ can be a subset of $\Pi$, while we require any three directions are not coplanar in subsection \ref{subsec-Unique-ManyDipoles}.
	
% lemma-upper bound of $f(z)$ - dipoles are in a particular plane
\begin{lemma}
\label{lemma-UpperBound-f(z)-InAPlane}
Given $M$ points $z_1, z_2, \cdots, z_M$ in $\Pi$ and observation direction set $\Theta_L\subset \Pi$.  Let $\Pi_{l,m}, 1\leq l\leq L, 1\leq m\leq M$ be the planes defined in \eqref{Definition-Planes}. For any point $z\in \Pi$ define by
\ben
g(z):= {\mathrm Card} \{\Pi_{l,m}\,|\, z\in \Pi_{l,m},\, 1\leq l\leq L,\, 1\leq m\leq M\}
\enn
the number of the planes $\Pi_{l,m}, 1\leq l\leq L, 1\leq m\leq M$ passing through $z$.	If any two observation directions in $\Theta_L$ are linearly independent, then
\ben
g(z)
\left\{
  \begin{array}{ll}
    \leq M, & \hbox{$z\in \mathbb{R}^3 \setminus \{z_1,z_2, \cdots, z_M\}$;} \\
    = L, & \hbox{$z\in \{z_1,z_2, \cdots, z_M\}$.}
  \end{array}
\right.
\enn
\end{lemma}
% begin of proof
\begin{proof}
It is obvious by definition that $g(z) =L$ for $z\in \{z_1,z_2, \cdots, z_M\}$.
	
For $z^* \in \Pi \setminus \{z_1,z_2, \cdots, z_M\}$, we show $g(z^*)\leq M$. Otherwise, there exists a point $z_m$ such that at least two planes $\Pi_{l_1,m}, \Pi_{l_2,m}$ pass through $z^*$, i.e.,
\ben
\hat{x}_{l_j} \cdot (z^* - z_m) = 0, \qquad j =1,2.
\enn
Noting that $\hat{x}_{l_1},\hat{x}_{l_2}, z^*, z_m\in \Pi$ and $\hat{x}_{l_1},\hat{x}_{l_2}$ are two linearly independent directions, we deduce from the above equality that $z^* =z_m $. This is a contradiction to $z^* \in \Pi \setminus \{z_1,z_2, \cdots, z_M\}$.

\end{proof}
% end of proof
%
%	Compare Lemma \ref{lemma-UpperBound-f(z)-InAPlane} with Lemma \ref{lemma-UpperBound-f(z)}, the upper bound of $f(z)$, where $z$ is not the location of dipoles, is smaller, which indicates less observation directions is needed.

Using Lemma \ref{lemma-UpperBound-f(z)-InAPlane} and following the proof of Theorem \ref{Theorem-ElecMagUnique-multi-frequency}, we immediately have the following unique result with less data.
	
% theorem-Uniqueness ManyDipoles InAPlane
\begin{theorem}
\label{Theorem-ElecMagInAPlane}
We consider $M$ dipoles located in a plane $\Pi$. Assume that the first $M_1$ dipoles are magnetic dipoles. If any two directions in $\Theta_L$ are linearly independent and $L>\max\{2M_1, 2(M-M_1)\}$, then
\ben
\{(z_m,q_m)\,|\, 1\leq m \leq M_1\}\quad\mbox{and}\quad \{(z_m,q_m)\,|\, M_1< m \leq M\}
\enn
are uniquely determined by multi-frequency electric far field patterns $E^\infty(\pm \hat{x},k), \hat{x}\in \Theta_L, k\in (0,\infty)$.
\end{theorem}
\section{Numerical algorithms}
\label{sec-Algorithm}
Following the idea in the proof of uniqueness results in the previous section, we introduce some numerical algorithms for locating the dipoles and reconstructing the corresponding polarization strengths.
%
%	In this section, we raise numerical methods that identify locations and polarization strengths of magnetic dipoles first. Then the electric dipoles is considered. In the end, we introduce an iterative algorithm that idetifis both electric and megnetic dipoles.
	
% subsection-identify-magnetic dipoles
\subsection{Magnetic dipoles}\label{subsec-Numerical-MagDipoles}

As the discussion in the uniqueness analyses, we begin with identifying the magnetic dipoles.
% subsubsection $z_m$
\subsubsection{Indicator for locating the magnetic dipoles}
Inspired by the proof of Theorem \ref{Theorem-ElecMagUnique-multi-frequency}, we define
\be
\label{eq-indicator-mag}
I_{mag}(z,K,\epsilon) := \frac{1}{L} \sum_{\hat{x}\in \Theta_L} T_{\epsilon} {(}\left|
F_{mag}(z,\hat{x},K)
\right|{)}
\en
with
\be
\label{eq-cut-function}
T_{\epsilon}(t) :=
\left\{
\begin{aligned}
&1, & t>\epsilon; \\
&0, & t\leq \epsilon
\end{aligned}
\right.
\en
and
\be
\label{eq-IntergralOfMag}
F_{mag}(z,\hat{x},K) := \frac{2 \pi}{K} \int_0^K \frac{1}{ik}\left[
e^{ik\hat{x}\cdot z} E^\infty(\hat{x},k) - e^{-ik\hat{x}\cdot z} E^\infty(-\hat{x},k)
\right]{\rm d}{k},\quad \hat{x}\in \Theta_L, K>0.
\en
By the analysis in the proof of Theorem \ref{Theorem-ElecMagUnique-multi-frequency}, $I_{mag}(z,K,\epsilon)$ goes to zero outside the planes $\Pi_{l,M}, 1\leq l\leq L, 1\leq m\leq M_1$ when $K$ tends to infinity. Furthermore, with the increase of the number $L$ of the properly chosen observation directions, $I_{mag}(z,K,\epsilon)$ tends to one for $z\in\{z_1, z_2, \cdots, z_{M_1}\}$. Therefore, the indicator \eqref{eq-indicator-mag} can be used to locate the magnetic dipoles with sufficient many properly chosen observation directions and sufficiently large $K$.

Finally, we want to give a remark on the auxiliary function \eqref{eq-cut-function}. From the point of view of numerical computations, the strength $\tau_m$ plays the role of weight for the $m-$th dipole, therefore the dipoles with quite small strengths are difficult to be located.
The introduction of the auxiliary function \eqref{eq-cut-function} is used to balance these weights and to avoid missing the dipoles with small strengths.
Practically, we may choose the cut-off value $\epsilon$ such that $|q_m|>\epsilon$ for all $m=1,2,\cdots, M_1$.

\subsubsection{Identifying the polarization strengths of magnetic dipoles}

Having located all the magnetic dipoles, we now consider the corresponding polarization strengths. For a magnetic dipole $(z_{m^*}, q_{m^*})$, we take two linearly independent directions $\hat{x}_{m^*},\hat{y}_{m^*}\in\Theta_L$ such that
\be\label{magdirections}
\hat{x}_{m^*}\cdot (z_m - \cdot z_{m^*})\neq 0,\quad \hat{y}_{m^*}\cdot (z_m - \cdot z_{m^*})\neq 0, \quad\forall\, m\neq {m^*}, 1\leq m\leq M_1.
\en
In terms of \eqref{qmmagformula}, we have
\be\label{qmmagast}
q_{{m^*}}
&=& \frac{\hat{y}_{m^*} \times \hat{x}_{m^*}}{|\hat{y}_{m^*}\times \hat{x}_{m^*}|^2} \cdot \Big(
F_{mag} (z_{m^*}, \hat{y}_{m^*})+\hat{y}_{m^*} \times [ \hat{x}_{m^*} \times F_{mag} (z_{m^*}, \hat{x}_{m^*}) ]\Big) \hat{x}_{m^*}\cr
&& -\hat{x}_{m^*} \times F_{mag} (z_{m^*}, \hat{x}_{m^*}), \quad 1\leq m^*\leq M_1.
\en
Numerically, we replace $F_{mag}(z_{m^*}, \hat{x}), F_{mag}(z_{m^*}, \hat{y})$ by $F_{mag}(z_{m^*}, \hat{x}, K), F_{mag}(z_{m^*}, \hat{y}, K)$ and set
\be\label{qmmagastK}
q_{{m^*} ,K}
&:=& \frac{\hat{y}_{m^*} \times \hat{x}_{m^*}}{|\hat{y}_{m^*}\times \hat{x}_{m^*}|^2} \cdot \Big(
F_{mag} (z_{m^*}, \hat{y}_{m^*},K)+\hat{y}_{m^*} \times [ \hat{x}_{m^*} \times F_{mag} (z_{m^*}, \hat{x}_{m^*},K) ]\Big) \hat{x}_{m^*}\cr
&& -\hat{x}_{m^*} \times F_{mag} (z_{m^*}, \hat{x}_{m^*},K), \quad 1\leq m^*\leq M_1.
\en

\begin{theorem}
\ben
\left| q_{m^*} - q_{{m^*},K} \right| = O\left(
\frac{1}{K}
\right),\quad K\rightarrow\infty,\quad 1\leq m^*\leq M_1.
\enn
\end{theorem}
\begin{proof}
With the help of \eqref{eq-integralOfMag} we have
\ben
\label{eq-ErrorAna_FMag}
&&\left|F_{mag}(z_{m^*},\hat{x}, K) - \hat{x} \times q_{m^*}\right| \cr
&=& \left|\sum_{m\neq m^*, m=1}^{M_1}  \int_{-K}^K\frac{e^{ik \hat{x}\cdot (z_{m^*}-z_m)}}{2K}
 {\rm d}{k}\  \hat{x} \times q_m  + \sum_{m=M_1+1}^M  \int_{0}^K \frac{i\sin(k \hat{x}\cdot(z_{m^*}-z_m))}{K}
{\rm d}{k}\  \hat{x} \times (q_m \times \hat{x})
\right| \cr
&\leq& \frac{1}{K} \sum_{m\neq m^*, m=1}^M \frac{|q_m|}{|\hat{x}\cdot (z_{m^*}-z_m)|}.
\enn
From this, comparing \eqref{qmmagast} and \eqref{qmmagastK} we have
\ben\label{Error-q-mag}
\left| q_{m^*} - q_{{m^*},K} \right|
&\leq& |F_{mag}(z_{m^*}, \hat{x},K) - \hat{x}\times q_{m^*}|\cr
&& +\frac{1}{|\hat{y}\times \hat{x}|} \Big(|F_{mag}(z_{m^*}, \hat{x},K) - \hat{x}\times q_{m^*}|+|F_{mag}(z_{m^*}, \hat{y},K) - \hat{y}\times q_{m^*}|\Big) \cr
&\leq& \frac{1+|\hat{y}\times \hat{x}|}{K|\hat{y}\times \hat{x}|} \sum_{m\neq m^*, m=1}^M \frac{|q_m|}{|\hat{x}\cdot (z_{m^*}-z_m)|}  +
\frac{1}{K|\hat{y}\times \hat{x}|} \sum_{m\neq m^*, m=1}^M \frac{|q_m|}{|\hat{y}\cdot (z_{m^*}-z_m)|}.\quad
\enn
This completes the proof.
\end{proof}

Letting $K\rightarrow\infty$ we see that
\ben
\lim_{K\to \infty} q_{{m^*},K} = q_{m^*}, \quad 1\leq m^*\leq M_1.
\enn
Therefore, $ q_{{m^*},K}$ can be an approximation to $q_{m^*}$ for large $K$.

%	when $\hat{x},\hat{y}$ is fixed.

% subsection-identify-electric dipoles
\subsection{Electric dipoles}\label{subsec-Numerical-ElecDipoles}
Having identified the magnetic dipoles, we look for the electric dipoles. We omit the analysis since it's similar to the case of magnetic dipoles.
	
\subsubsection{Indicator for the locations of electric dipoles}
We define
\be
\label{eq-indicator-elec}
I_{elec}(z,K,\epsilon) := \frac{1}{L} \sum_{\hat{x}\in \Theta_L} T_{\epsilon} {(}\left|
F_{elec}(z,\hat{x},K)
\right|{)}
\en
	as the indicator function for locating the electric dipoles, where
$T_\epsilon$ defined in \eqref{eq-cut-function}
and
\be
\label{eq-IntergralOfelec}
F_{elec}(z,\hat{x},K) := \frac{2 \pi}{K} \int_0^K \frac{1}{ik}\left[
e^{ik\hat{x}\cdot z} E^\infty(\hat{x},k) + e^{-ik\hat{x}\cdot z} E^\infty(-\hat{x},k)
\right]{\rm d}{k},\quad \hat{x}\in \Theta_L, K>0.
\en

Moreover, we define
\be
\label{eq-limit-Felec}
F_{elec} (z, \hat{x}) := \lim_{K\to \infty} F_{elec} (z, \hat{x},K).
\en

% subsubsection $q_m$
\subsubsection{Identifying the polarization strengths of electric dipoles}
Given all the location of electric dipoles $\{z_{M_1+1}, z_{M_1+2}, \cdots, z_M\}$, for the ${m^*}$th electric dipole, we take two linearly independent directions $\hat{x}_{m^*}, \hat{y}_{m^*}$ such that
\be\label{elecdirections}
\hat{x}_{m^*}\cdot (z_m - \cdot z_{m^*})\neq 0,\quad \hat{y}_{m^*}\cdot (z_m - \cdot z_{m^*})\neq 0, \quad\forall\, m\neq {m^*}, M_1< m\leq M.
\en
In view of \eqref{qmelecformula}, we have
\be\label{qmelecast}
q_{m^*}
&=& \frac{\hat{y}_{m^*} \times \hat{x}_{m^*}}{|\hat{y}_{m^*}\times \hat{x}_{m^*}|^2} \cdot \Big(\hat{y}_{m^*} \times F_{elec} (z_{m^*}, \hat{y}_{m^*})-\hat{y}_{m^*} \times F_{elec} (z_{m^*}, \hat{x}_{m^*}) \Big)\hat{x}_{m^*}\cr
&&+ F_{elec} (z_{m^*}, \hat{x}_{m^*}),\quad M_1< {m^*}\leq M.
\en
Replacing $F_{elec} (z_{m^*}, \hat{x}_{m^*}),F_{elec} (z_{m^*}, \hat{y}_{m^*})$ by $F_{elec} (z_{m^*}, \hat{x}_{m^*},K),F_{elec} (z_{m^*}, \hat{y}_{m^*},K)$, respectively, we get the approximation of $q_{m^*}$ by
\be\label{qmelecastK}
q_{m^* ,K}
&:=& \frac{\hat{y}_{m^*} \times \hat{x}_{m^*}}{|\hat{y}_{m^*}\times \hat{x}_{m^*}|^2} \cdot \Big(\hat{y}_{m^*} \times F_{elec} (z_{m^*}, \hat{y}_{m^*},K)-\hat{y}_{m^*} \times F_{elec} (z_{m^*}, \hat{x}_{m^*},K) \Big)\hat{x}_{m^*}\cr
&&+ F_{elec} (z_{m^*}, \hat{x}_{m^*},K),\quad M_1< {m^*}\leq M.
\en
Similarly, we have $\left| q_{m^*} - q_{{m^*},K} \right| = O\left(\frac{1}{K}\right)$ as $K\rightarrow\infty$.
% subsection-algorithm
\subsection{Algorithm for all dipoles}
\label{subsec-Algorithm}

We combine all the indicators and formulas introduced in the previous subsections to form the following imaging algorithm.

\noindent\textbf{Imaging Algorithm for multiple mixed type dipoles.}
\
\begin{itemize}
\item Collect the multi-frequency sparse electric far field patterns $E^\infty(\pm \hat{x},k), \hat{x}\in \Theta_L, k\in (0,K)$.
\item Select a sampling region in $\mathbb{R}^3$ with a fine mesh containing all the dipoles.
\item Locate all the magnetic dipoles by plotting the indicator $I_{mag}(z,K,\epsilon)$ given in \eqref{eq-indicator-mag} with properly chosen $\epsilon>0$. Reconstruct the corresponding polarization strengths by the formula \eqref{qmmagastK}.
\item Locate all the electric dipoles by plotting the indicator $I_{elec}(z,K,\epsilon)$ given in \eqref{eq-indicator-elec} with properly chosen $\epsilon>0$. Reconstruct the corresponding polarization strengths by the formula \eqref{qmelecastK}.
\end{itemize}

% ----------------------------------------------------------
% -------------------- New Section -------------------------
% Numerical Example
\section{Numerical example}
\label{sec-NumericalExample}

In this section, we present some numerical simulations to verify the effectiveness and robustness of the proposed numerical algorithm. Fix the wave number $k$, denote the electric far field pattern at wave number $k$ by $F(k) := \left( E^\infty(\hat{x}_l, k)\right)_{ 1\leq l\leq L} \in \mathbb{C}^{L\times 3}$. We perturb $F(k)$ by random noise using
\be
\label{Definition-FarP-withNoise}
F^\delta (k) = F(k) + \delta | F(k)| \frac{R_1 + R_2 i}{|R_1 + R_2 i|},
\en
where $R_1, R_2$ are two $L\times 3$ matrixes containing values drawn from a normal distribution with mean zero and standard derivation one. The
value of $\delta$ used in our code is $\delta = | F^\delta(k) - F(k)| / |F(k)|$ and so presents the relative error.

In the simulations, we used a grid $\mathcal{G}$ of $M\times M \times M$ equally spaced sampling points on some rectangle $[-c,c]\times[-c,c]\times[-c,c]$.
For each point $z \in \mathcal{G}$, we locate the dipoles by indicators $I^{\rho}_{mag}(z,K,\epsilon)$ and $I^{\rho}_{elec}(z,K,\epsilon)$ with $K=100$ and $\epsilon=0.2$.
Here, $\rho\geq1$ is an artificial selected parameter to enhance the resolution.

If not otherwise stated, the observation directions $\hat{x}_l=(x^1_l, x^2_l, x^3_l)$ in $\Theta_L$ are selected by Fibonacci lattices:
\be
\label{Definition-evenly-direction}
x^3_l  := 1-\frac{2l}{L};	\quad
x^1_l  := \sqrt{1-(x^3_l)^2} \cos(2\pi l \phi);	\quad
x^1_2  := \sqrt{1-(x^3_l)^2} \sin(2\pi l \phi),\quad 1\leq l\leq L.
\en
Here, $\phi = \frac{\sqrt{5}-1}{2}$ is the golden ratio to ensure that the lattices are evenly spaced \cite{Alvaro}. Note that we consider the Fibonacci lattices just because we guess any three Fibonacci lattices are not coplanar. However, to our best knowledge, this is still not be proved.

% subsection-Identify location and strgenth vector of dipoles
\subsection{Identifying the locations and the polarization strengths of the dipoles}
\label{subsec-NumericalExampla-manyMag}

\begin{figure}[h!]
    \centering
    \begin{tabular}{cc}
%  x=0
        \subfigure[$z^1=0,  I^4_{mag}(z)$]{
            \label{fig-EAM-Mag-x=0}
            \includegraphics[width=.50\textwidth]{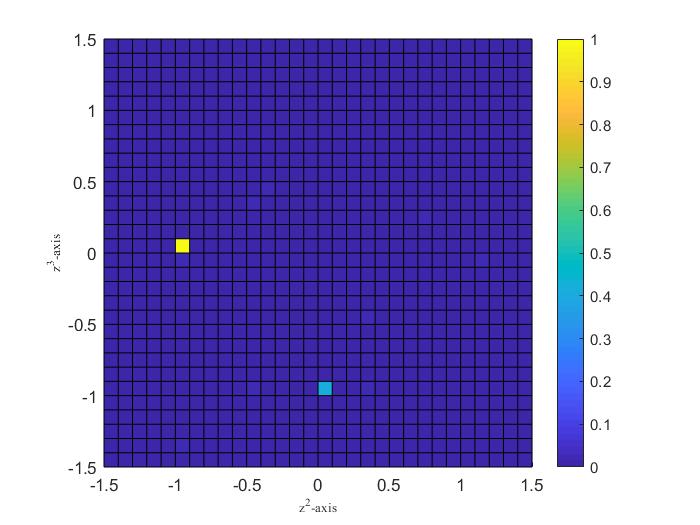}
        } \hspace{0em} &
        \subfigure[$z^1=0,   I^4_{elec}(z)$]{
            \label{fig-EAM-Elec-x=0}
            \includegraphics[width=.50\textwidth]{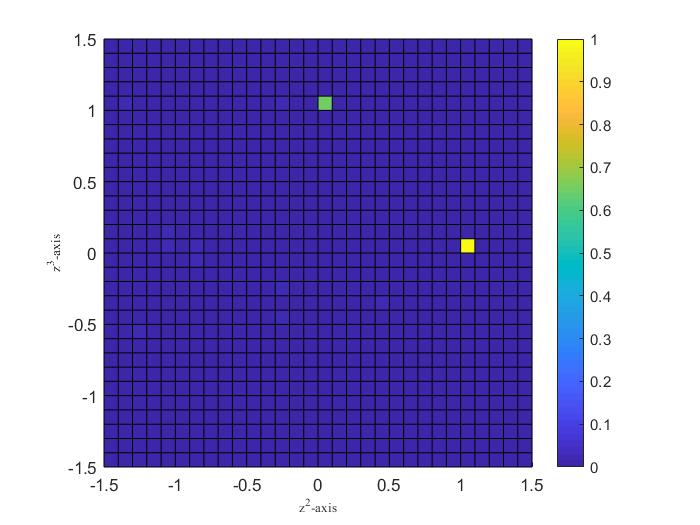}
        } \\
%  y=0
        \subfigure[$z^2=0,   I^4_{mag}(z)$]{
            \label{fig-EAM-Mag-y=0}
            \includegraphics[width=.50\textwidth]{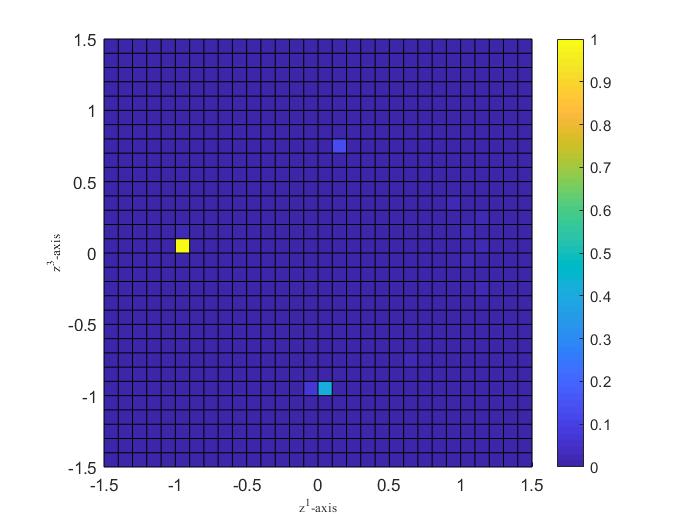}
        } \hspace{0em} &
        \subfigure[$z^2=0,   I^4_{elec}(z)$]{
            \label{fig-EAM-Elec-y=0}
            \includegraphics[width=.50\textwidth]{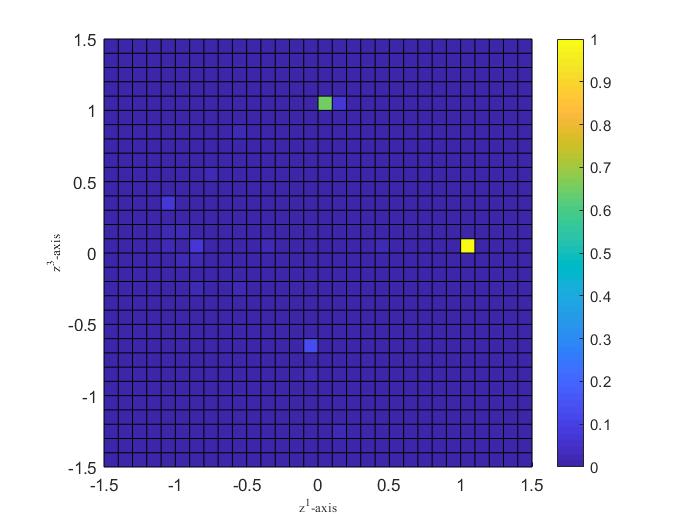}
        } \\
%  z=0
        \subfigure[$z^3=0,   I^4_{mag}(z)$]{
            \label{fig-EAM-Mag-z=0}
            \includegraphics[width=.50\textwidth]{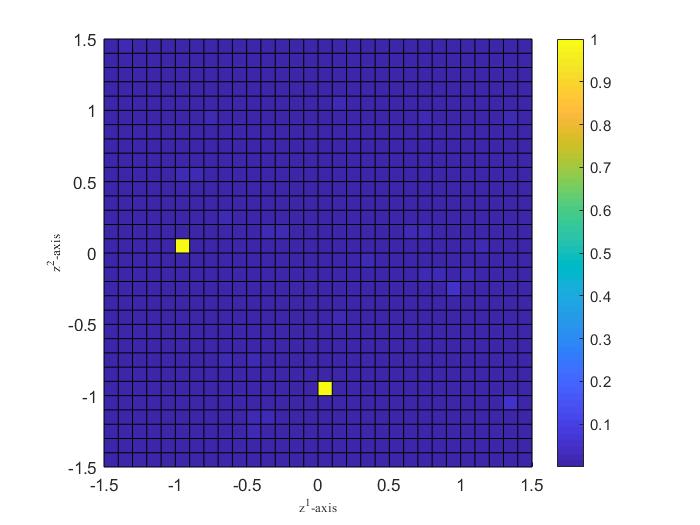}
        } \hspace{0em} &
        \subfigure[$z^3=0,   I^4_{elec}(z)$]{
            \label{fig-EAM-Elec-z=0}
            \includegraphics[width=.50\textwidth]{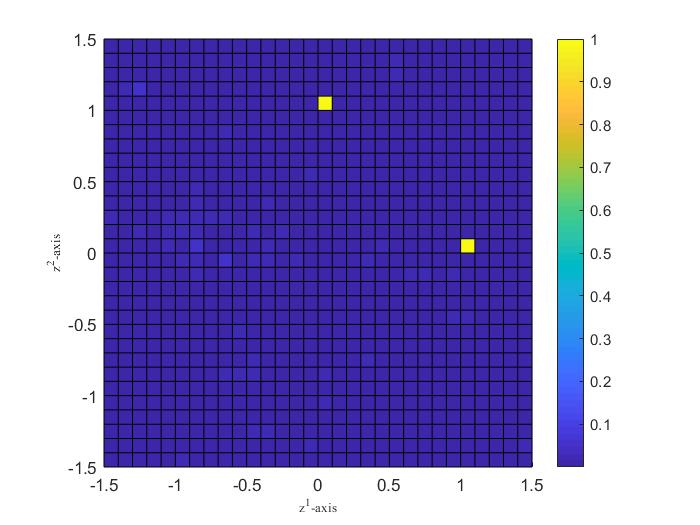}
        } \\
    \end{tabular}
    \caption{Location reconstructions in the coordinate planes with $L=10$ and $10\%$ relative noise.}
    \label{fig-MagAndElec}
\end{figure}

In the first example, we consider the mixed type point sources with $3$ magnetic dipoles and $3$ electric dipoles. The locations and polarization strengths are shown in the second and third columns of Table \ref{table-Retrival-ElecMag-10noise}.
To locate the dipoles, we consider the research domain $[-1.5,1.5]\times[-1.5,1.5]\times[-1.5,1.5]$ with $31\times 31\times 31$ equally spaced sampling points.
Figure \ref{fig-MagAndElec} shows the location reconstructions in different hyper-planes $z^i = 0, i=1,2,3$.
As shown in Table \ref{table-Retrival-ElecMag-10noise}, in the hyperplane $z^1=0$, there are two magnetic dipoles located at $(0, -1, 0)$ and $(0, 0, -1)$, respectively. Obviously, Figure \ref{fig-MagAndElec}(a) shows that these two magnetic dipoles are well captured by the indicator $I^{4}_{mag}$. Meanwhile, as shown in Figure \ref{fig-MagAndElec}(b), the two electric dipoles located in the hyperplane $z^1=0$ are clearly reconstructed by the indicator $I^{4}_{elec}$. Considering $10\%$ relative noise in the measurement data, the location reconstructions are quite stable. This can also be seen in Figure \ref{fig-MagAndElec}(c-f) for the other
reconstructions in the coordinate planes.

Having located all the dipoles, we compute the corresponding polarization strengths by the formulas \eqref{qmmagastK} and \eqref{qmelecastK}.
For each strength reconstruction, we have used the electric far field patterns with $K=200$ at two selected observation directions satisfying \eqref{magdirections} or \eqref{elecdirections}.
The fourth column of Table \ref{table-Retrival-ElecMag-10noise} shows the reconstructed polarization strengths. Considering $10\%$ relative noise in the measurements again, the relative errors for the strengths shown in the fifth column of Table \ref{table-Retrival-ElecMag-10noise} are acceptable.

\begin{table}
     \centering
     \begin{tabular}{|c|c|c|c|c|}
     \hline \multicolumn{1}{|c|}{\textbf{Type}} &
     \multicolumn{1}{c|}{\textbf{$z_m$}} &
     \multicolumn{1}{c|}{\textbf{True $q_m$}} &
     \multicolumn{1}{c|}{\textbf{Reconstructed $q_{m}$}} &
     \multicolumn{1}{c|}{\textbf{RE}} \\ \hline
  1st Magnetic  & (-1, 0, 0)   &  (1, 1, -1) & (1.02+0.00i, 1.00+0.02i, -1.01-0.04i) & 2.89\%\\
 \hline
  2nd Magnetic  & (0, -1, 0)   &  (-0.5, 0, 1) & (-0.49+0.03i, -0.01-0.01i, 1.00-0.02i) & 3.58\%\\
 \hline
  3rd Magnetic  & (0, 0, -1)   &  (-1, 0.2, 0) & (-1.00+0.01i, 0.15+0.01i, 0.01-0.02i) & 5.55\%\\
 \hline
  1st Electric  & (1, 0, 0)   &  (1, 1, 1) & (0.97-0.04i, 0.99+0.04i, 1.01-0.01i) & 3.83\%\\
 \hline
  2nd Electric  & (0, 1, 0)   &  (0.5, 0, 1) & (0.49-0.01i, 0.99+0.04i, 1.01-0.01i) & 1.55\%\\
 \hline
  3rd Electric  & (0, 0, 1)   &  (1, 0.2, 0) & (0.99+0.00i, 0.16+0.04i, 0.00+0.03i) & 6.35\%\\
 \hline
     \end{tabular}
     \caption{The first example with 3 magnetic dipoles and 3 electric dipoles. In the fourth column, we present the reconstructed $q_m$ with $10\%$ relative noise in the measurements. The fifth column shows the relative error (RE) of the reconstructed polarization strengths $q_m$.}
  \label{table-Retrival-ElecMag-10noise}
\end{table}

% Graph - Values of Indicatoe - Many Magnetic Dipoles -Directions in A Plane

% subsection-Identify location of Many Magnetic dipoles
\subsection{Locating multiple magnetic dipoles using different observation directions}
\label{subsec-NumericalExampla-manyMag}
In the second example, we consider locating  $M_1=19$ magnetic dipoles. The positions and polarization strengths vector are presented in Table \ref{table-ManyMagDipoles-Data}. All the dipoles are located in the hyper-plane $\Pi := \{(z^1, z^2, z^3)\in \mathbb{R}^3| z^3=0\}$. We consider the following three
observation direction sets:
\begin{itemize}
  \item A single pair of directions $\pm\hx$.
  \item $L$ directions in the unit sphere $S^2$ chosen by the Fibonacci lattices \eqref{Definition-evenly-direction};
  \item $L$ equally distributed directions in the hyper-plane $\Pi$, i.e.,
\be\label{Definition-evenly-direction-InAPlane}
x^1_l  = \cos(\pi l/L),	\quad x^1_l  =  \sin(\pi l/L),\quad x^3_l = 0,\qquad 1\leq l\leq L.
\en
\end{itemize}

Due to the a priori information that all the dipoles are located in $\Pi$, we consider the research domain $[-2, 2]\times[-2, 2]\times 0$ with $41\times 41$ equally spaced sampling points.

\begin{figure}[h!]
    \centering
    \begin{tabular}{cc}
        \subfigure[$\hat{x} = \pm (1,0,0)$]{
            \label{fig-SD-up}
            \includegraphics[width=.50\textwidth]{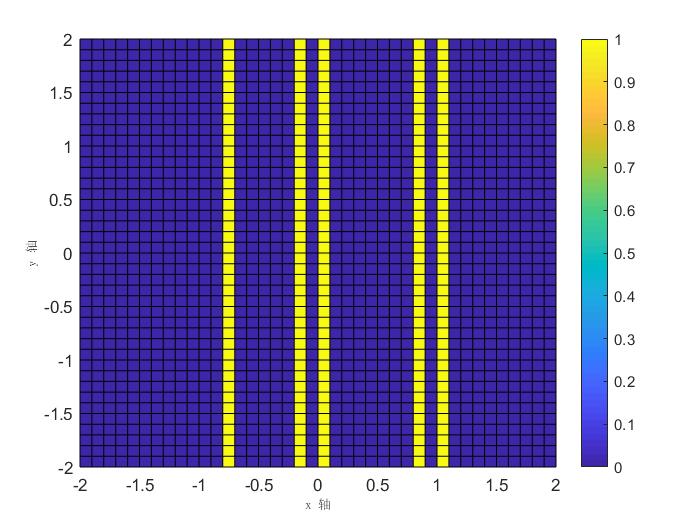}
        } \hspace{0em} &
        \subfigure[$\hat{x} = \pm (0,1,0)$]{
            \label{fig-SD-left}
            \includegraphics[width=.50\textwidth]{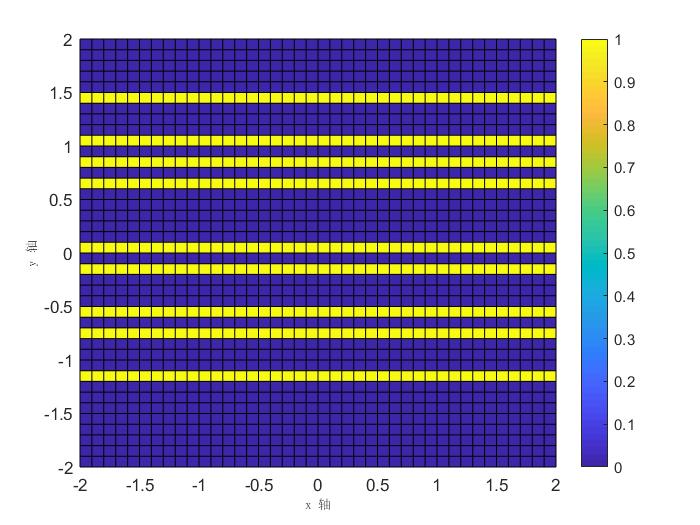}
        } \\
        \subfigure[$\hat{x} = \pm (\sqrt{2}/2, \sqrt{2}/2,0)$]{
            \label{fig-SD-UpLeft}
            \includegraphics[width=.50\textwidth]{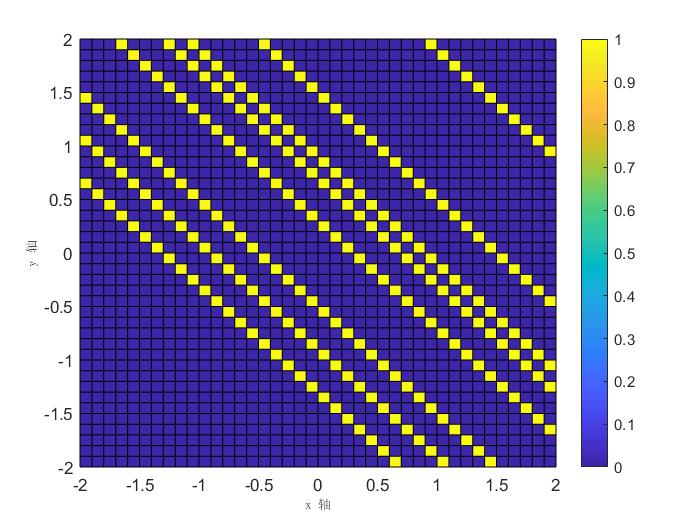}
        } \hspace{0em} &
        \subfigure[$\hat{x} = \pm (\sqrt{2}/2, -\sqrt{2}/2,0)$]{
            \label{fig-SD-UpRight}
            \includegraphics[width=.50\textwidth]{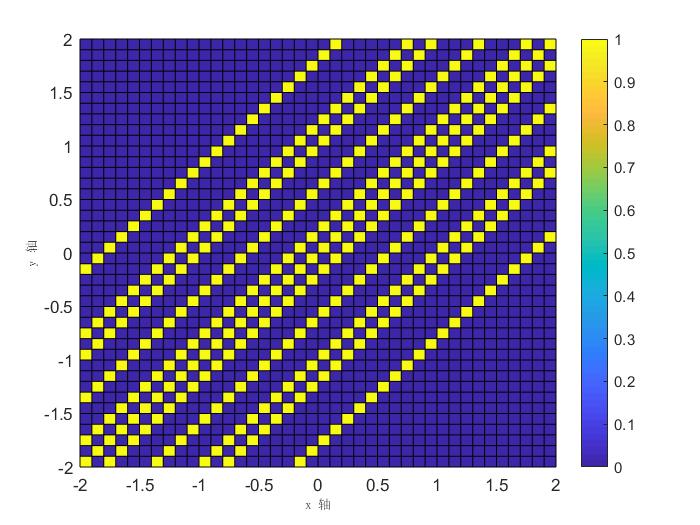}
        }\\
    \end{tabular}
    \caption{Reconstructions with a single pair of directions. $10\%$ noise is considered.}
    \label{fig-SingleDirection}
\end{figure}

Figure \ref{fig-SingleDirection} shows the reconstructions with a single pair of directions $\pm\hx$. The highlighted lines show that there must be at least one dipole located in the lines. Conversely, the line passing through the dipoles may be missed. As shown in the Figure \ref{fig-SingleDirection}($a$), we have observe only five lines $z^1=-0.8, -0.2, 0, 0.8$ and $1$. The lines passing through the dipoles with polarization strengths paralleling to $\hx=(1,0,0)$ are missed. This is obvious by noting the representation \eqref{eq-FarFOfElecAndMagDipoles} of the electric far field pattern.

% Table-Data of Magnetic Dipoles
% Table-Data of Magnetic Dipoles
\begin{table}
     \centering
     \begin{tabular}{|c|c|c|c|} \hline
    \multicolumn{1}{|c|}{\textbf{Order}} &
    \multicolumn{1}{c|}{\textbf{Location $z_m$}} &
    \multicolumn{1}{c|}{\textbf{Polarization strength $q_m$}} \\ \hline
1  & (1.4, 1.4, 0)  & (1.20+1.49i, 0, 0)    \\ \hline
2  & (0.8, 1.4, 0)  & (0, -1.40-0.63i, 0)   \\ \hline
3  & (0.4, 1.0, 0)  & (1.00-0.96i, 0, 0)    \\ \hline
4  & (-0.2, 1.0, 0) & (1.20+0.84i, 0, 0)    \\ \hline
5  & (-0.8, 1.0, 0) & (0, 0.83-1.41i, 0)    \\ \hline
6  & (-1.2, 0.6, 0) & (0.90+1.43i, 0, 0)    \\ \hline
7  & (0, 0.6, 0)   & (0, 0, 1.35-0.97i)   \\ \hline
8  & (-0.8, -0.2, 0) & (0, 0, 1.19+1.28i)    \\ \hline
9  & (-1.2, -0.2, 0) & (1.45-0.58i, 0, 0)    \\ \hline
10  & (-0.8, -0.6, 0) & (0, 0, 0.67+1.44i)    \\ \hline
11  & (-0.6, -0.8, 0) & (-1.08-1.47i, 0, 0)   \\ \hline
12  & (-0.2, -1.2, 0) & (1.00+1.39i, 0, 0)    \\ \hline
13  & (-0.2, -0.8, 0) & (0, -1.16+1.44i, 0)    \\ \hline
14  & (0.6, 0, 0)   & (-0.52-0.70i, 0, 0)   \\ \hline
15  & (0.6, -1.2, 0) & (1.02+1.15i, 0, 0)    \\ \hline
16  & (1.0, -0.8, 0) & (0.77+0.79i, 0, 0)    \\ \hline
17  & (1.0, -0.2, 0) & (1.42-1.07i, 0, 0)    \\ \hline
18  & (1.0, 0.4, 0)  & (0, 0.82+0.87i, 0)    \\ \hline
19  & (1.4, 0.8, 0)  & (-1.32-0.65i, 0, 0)    \\ \hline
     \end{tabular}
    \caption{The 19 magnetic dipoles to be located in the second example. }
     \label{table-ManyMagDipoles-Data}
\end{table}

% Graph-Value of Indicator-Just one direction
	
\begin{figure}[h!]
    \centering
    \begin{tabular}{cc}
%  L=10
        \subfigure[$L=10$, $10\%$ noise]{
            \label{fig-MD-10dire-10noise}
            \includegraphics[width=.50\textwidth]{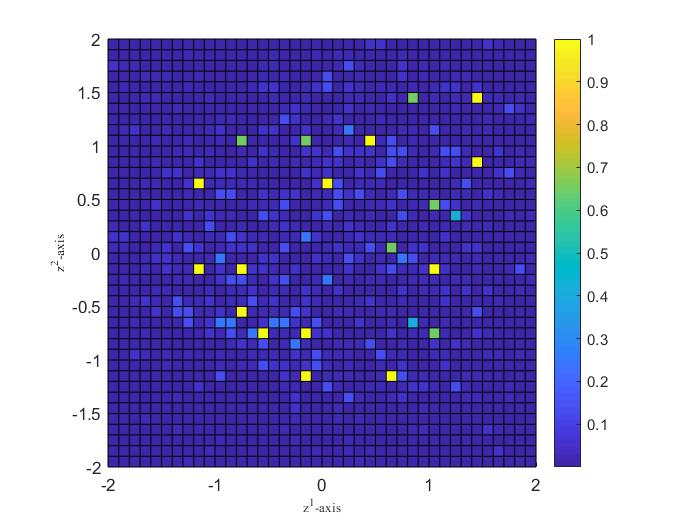}
        } \hspace{0em} &
        \subfigure[$L=10$, $20\%$ noise]{
            \label{fig-MD-10dire-20noise}
            \includegraphics[width=.50\textwidth]{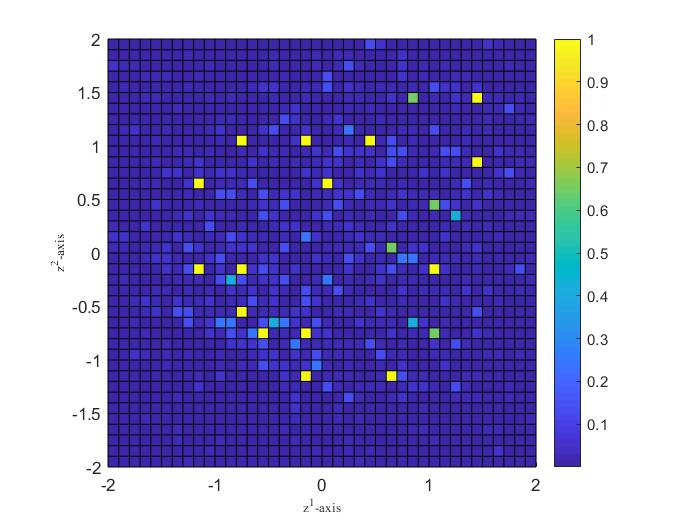}
        } \\
%  L=20
        \subfigure[$L=20$, $10\%$ noise]{
            \label{fig-MD-20dire-10noise}
            \includegraphics[width=.50\textwidth]{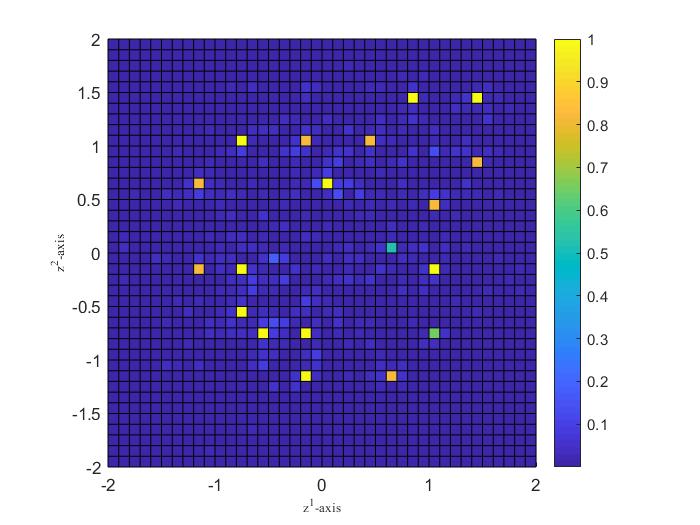}
        } \hspace{0em} &
        \subfigure[$L=20$, $20\%$ noise]{
            \label{fig-MD-20dire-20noise}
            \includegraphics[width=.50\textwidth]{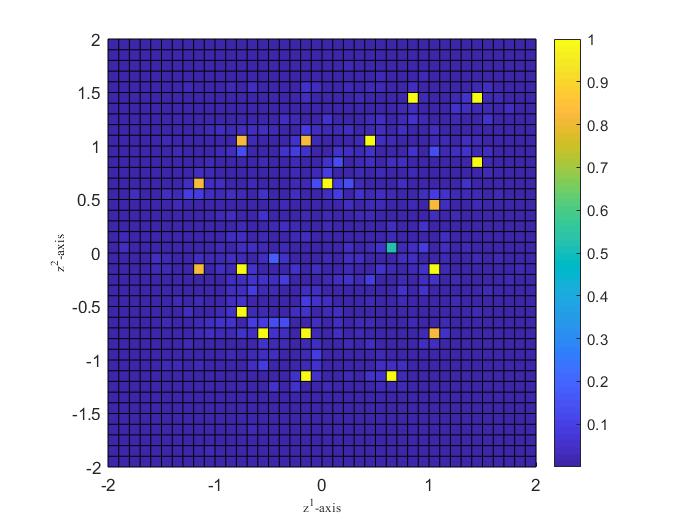}
        } \\
%  L=40
        \subfigure[$L=40$, $10\%$ noise]{
            \label{fig-MD-40dire-10noise}
            \includegraphics[width=.50\textwidth]{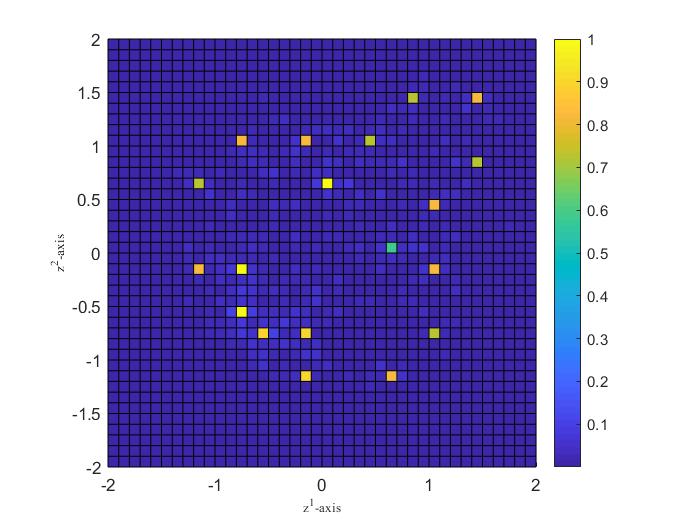}
        } \hspace{0em} &
        \subfigure[$L=40$, $20\%$ noise]{
            \label{fig-MD-40dire-20noise}
            \includegraphics[width=.50\textwidth]{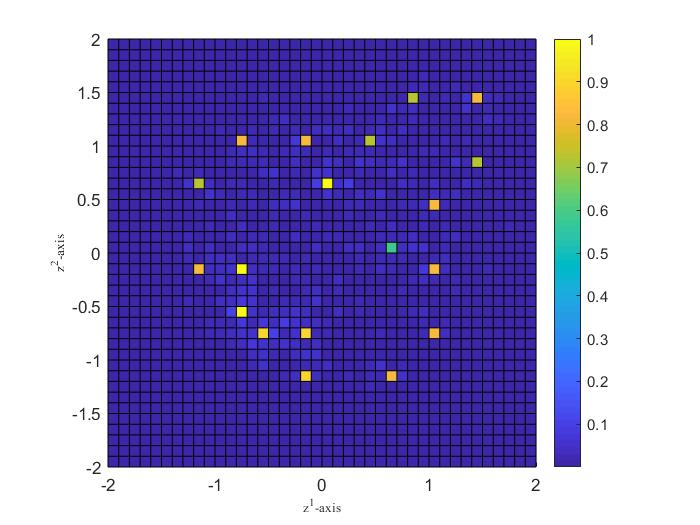}
        } \\
    \end{tabular}
    \caption{Reconstructions by the second direction set \eqref{Definition-evenly-direction} in $S^2$. }
    \label{fig-ManyDirection}
\end{figure}

Figure \ref{fig-ManyDirection} shows the results with more observation directions in $S^2$.
Obviously, with the increase of the number of the observation directions,
the unknown dipoles are clearly located  and the false positions are removed. Note that $L=40<4M_1=76$, i.e., we have used less data than those needed in the uniqueness Theorem \ref{Theorem-ElecMagUnique-multi-frequency}.

\begin{figure}[h!]
    \centering
    \begin{tabular}{cc}
%  L=10
        \subfigure[$L=10$, $10\%$ noise]{
            \label{fig-MDIAP-10dire-10noise}
            \includegraphics[width=.50\textwidth]{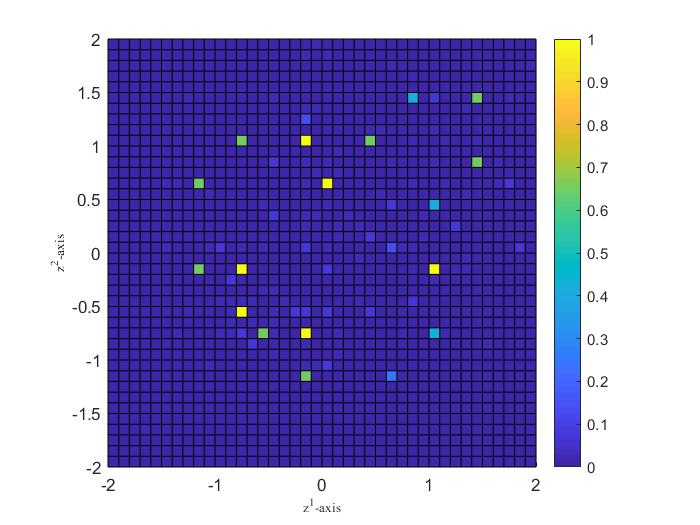}
        } \hspace{0em} &
        \subfigure[$L=10$, $20\%$ noise]{
            \label{fig-MDIAP-10dire-20noise}
            \includegraphics[width=.50\textwidth]{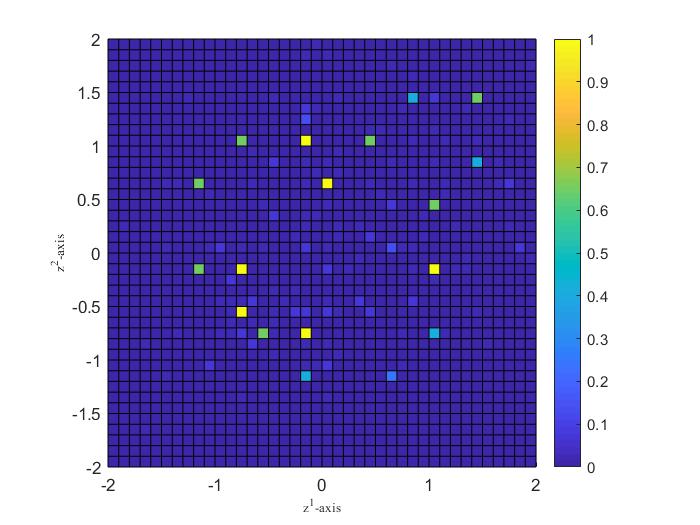}
        } \\
%  L=20
        \subfigure[$L=20$, $10\%$ noise]{
            \label{fig-MDIAP-20dire-10noise}
            \includegraphics[width=.50\textwidth]{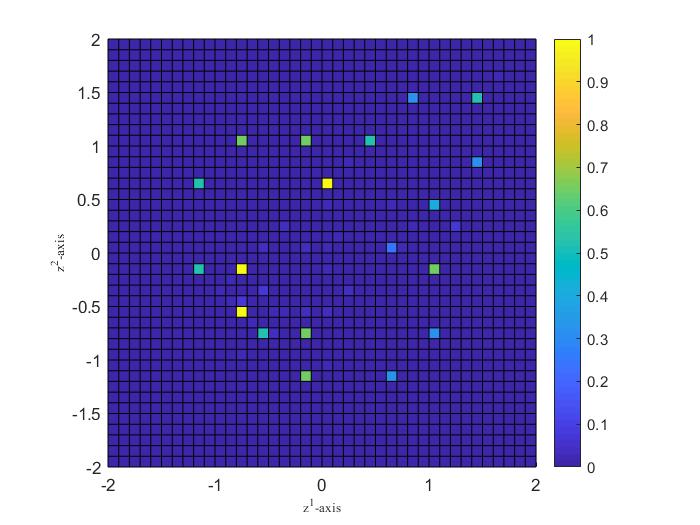}
        } \hspace{0em} &
        \subfigure[$L=20$, $20\%$ noise]{
            \label{fig-MDIAP-20dire-20noise}
            \includegraphics[width=.50\textwidth]{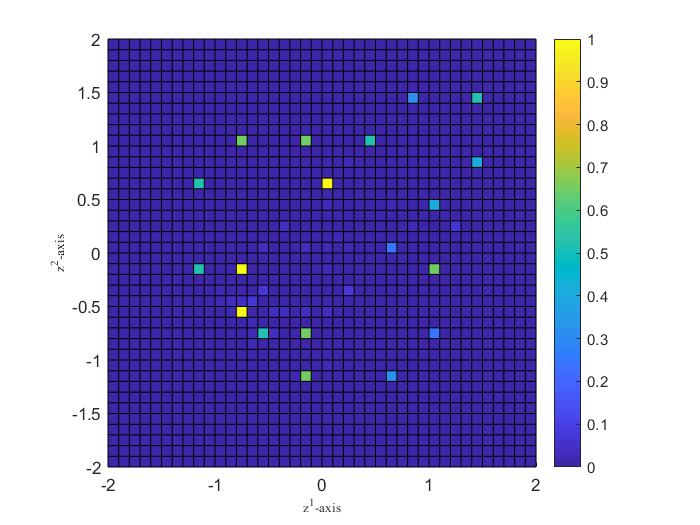}
        } \\
%  L=40
        \subfigure[$L=40$, $10\%$ noise]{
            \label{fig-MDIAP-40dire-10noise}
            \includegraphics[width=.50\textwidth]{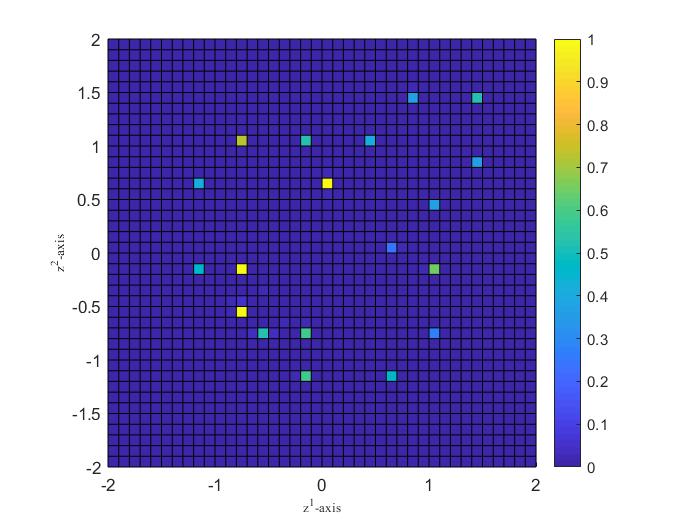}
        } \hspace{0em} &
        \subfigure[$L=40$, $20\%$ noise ]{
            \label{fig-MDIAP-40dire-20noise}
            \includegraphics[width=.50\textwidth]{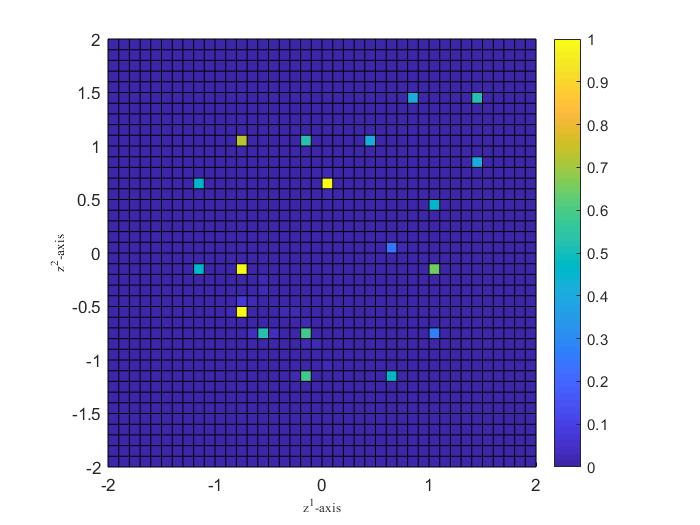}
        } \\
    \end{tabular}
    \caption{Reconstructions by the third direction set \eqref{Definition-evenly-direction-InAPlane} in $\Pi$}
    \label{fig-ManyDirection-InAPlane}
\end{figure}

Once we know that all the dipoles are in the plane $\Pi$, by Theorem \ref{Theorem-ElecMagInAPlane}, the locations of magnetic dipoles can be determined from electric far field pattern in the directions in $\Pi$. Figure \ref{fig-ManyDirection-InAPlane} shows the reconstruction using the third direction set.
Obviously, with the same direction number $L$, the reconstructions in Figure \ref{fig-ManyDirection-InAPlane} are better than those in Figure \ref{fig-ManyDirection}.

% Graph-Results of Indicator- Many Magnetic Dipols in a Plane

% Graph - Values of Indicatoe - Many Magnetic Dipoles -Directions in A Plane

% ----------------------------------------------------------
% -------------------- New Section -------------------------
% Conclude this article
\section{Conclusion and remark}
\label{sec-Conclution}
%We have proposed some novel techniques to identify the mixed type dipoles. Based on the uniqueness arguments, we have also introduced the corresponding indicator functions to locate the dipoles and proposed the formulas for computing their polarization strengths. The numerical examples further verify our theoretical analyses and numerical algorithm.
We have studied the uniqueness and numerical algorithm for identifying the mixed type dipoles. The novel geometrical arguments and ingenious integrals of the multi-frequency sparse electric far field patterns with properly chosen functions are key to the analyses. The numerical examples further verify our theoretical analyses and numerical algorithm. Precisely, the dipoles are well located and distinguished by the proposed indicator functions. Furthermore, the corresponding polarization strengths can also be reconstructed by the proposed formulas.

Similar ideas can also be applied to the multi-frequency electric fields taken at sparse sensors. However, more complex geometrical discussion will be involved. As mentioned in the introduction part, the dipoles can be viewed as the fundamental solution to the Maxwell equations. Thus we expect that our indicator functions are also applicable for determining the extended sources.

\section*{Acknowledgement}

The research of X. Liu is supported by the NNSF of China grant 11971471 and the Youth Innovation Promotion Association, CAS.

\bibliographystyle{SIAM}

\end{document}